\documentclass[12pt]{article}

\usepackage{amsmath,comment,lscape}
\usepackage{amsfonts, authblk}
\usepackage{amssymb}
\usepackage{graphicx}

\newtheorem{dfn}{Definition}[section]

\usepackage{caption}
\usepackage{subcaption}
\usepackage{xcolor}

\title{\textbf{A Fractional Model of Abalone Growth using Adomian Decomposition Method}}
\author[1,2]{Marliadi Susanto}
\author[2]{Nadihah Wahi}
\author[3]{Adem Kilicman}
\affil[1]{Department of Mathematics, Faculty of Mathematics and Natural Science, Universitas Mataram, 83125, Indonesia}
\affil[2]{Department of Mathematics and Statistics, Faculty of Science, Universiti Putra Malaysia, 43400 Serdang, Selangor, Malaysia}
\affil[3]{School of Mathematical Sciences, College of Computing, Informatics and Mathematics, Universiti Teknologi MARA, 40450 Shah Alam, Selangor, Malaysia}
\date{\today}

\begin{document}
\maketitle
\begin{abstract}
This study is a modification of the McKendrick equation into a growth model with fractional order to predict the abalone length growth. We have shown that the model is a special case of Taylor's series after it was analysed using Adomian decomposition method and Caputo fractional derivative. By simulating the series with some fractional orders, the results indicate that the greater the fractional order of the model, the series values generated are greater as well. Moreover, the series that is close to the real data is the one with a fractional order of $0.5$. Therefore, the growth model with a fractional order provides more accuracy than a classical integer order. 

\end{abstract}
\textit{\textbf{Keywords}: Modified Model, ADM, Taylor's series.\\
\textit{\textbf{Type}: Mathematical Modelling}}

\section{Introduction}
Abalone (Haliotis asinina) is one of the key ethno-fauna of West Nusa Tenggara Province and nationally as a marine commodity for export. The presence of this mollusk has played an important role in the coastal community's economy, not only for local consumption or sale in local markets but also for export to several countries in Asia, Europe, and the United States. Abalone harvesting in the wild has been excessive, leading to a drastic decline in its population, which could threaten the sustainability of the species. Therefore, the abalone cultivation is being widely developed on an industrial scale, but the limited availability of seed stock  and and the slow growth rate of abalone impose great challenges \cite{[1]}. In this study, the growth of abalone is described mathematically since mathematical modelling is an essential tools in understanding the phenomena. Thomas Robert Malthus in 1798 declared that population in the world would increase over time exponentially and exceed its resources, which was expressed mathematically as:
\begin{equation}
\frac{dP}{dt}=rP
\end{equation} 
where the initial condition $P(0)=P_0$. By solving the Eq.(1), we obtain
\begin{equation}
P(t)=P_0 e^{rt}
\end{equation}
with $P(t)$ representing the number of population at time $t$, $P_0$ is the initial number of population, and $r$ denotes the intrinsic growth rate. The model was applied to predict the number of population in Taraba state, and it showed more realistic results than the logistic model demonstrating a higher $R^2$ value \cite{[2]}. Exponential growth model also provide approximation the growth of carcinoma and melanoma properly \cite{[3]}. In addition, the exponential model was used to describe the dynamical cell growth \cite{[4]} and in estimating electricity demand in Cameroon \cite{[5]}.

Furthermore, in 1926, McKendrick introduced the growth model in partial differential equation with age structure in one population. Let $w(s,t)$ be the density of population with age $s$ at time $t$, that is
\begin{equation}
\frac{\partial w(s,t)}{\partial t}+\frac{\partial w(s,t)}{\partial s}=-\mu(s,t)w(s,t), 
\end{equation}
subject to the initial condition
\begin{equation}
\\w(s,0)=v_0(s), \text{  and} \text{  w(0,t)}=\int^{s_2}_{s_1}m(s,t)w(s,t)ds,
\end{equation}
where the value of $\mu(s,t)$ indicates the death process and $m(s,t)$ represents fertility, which are time-dependent. The model is well-known as McKendrick equation \cite{[6]}.  The model was developed by Gourley with involving non-linear effects, that is, the model was combined with competition model to study larval competition based on age-structured \cite{[7]}. In addition, the equation was expanded into McKendrick-Von Foerster equation with higher order numerical scheme for singular mortality cases \cite{[8]}.

On the other hand, the McKendrick equation has been utilized in various fields, including the determination of optimal hiring and retirement ages for employees \cite{[9]}, developing goodwill model for estimating the duration of a product's life cycle \cite{[10]}, and modeling the transmission of Mycobacterium tuberculosis in Japan, influenced by visitor numbers \cite{[11]}. Similarly, the model was used to modify SIR and SEIR model by involving age structure in each the compartment \cite{[12]}. The theory of McKendrick also was applied in ecology problems including zooplankton, insect, mosquito, and tick populations \cite{[13]}-\cite{[16]}. Hence, the McKendrick model has become a foundation of partial differential equations models. Moreover, Arora (2023) stated that the Eq. (3) is a pure growth equation by assuming the $\mu=0$. The equation was analysed using Adomian decomposition method (ADM) which was obtaining the exponential function as a solution \cite{[17]}.

Therefore, this study proposes a new fractional growth model as a modification of Eq.(3) to predict the body length of abalone. Here the new model is analysed by ADM to obtain a particular of Taylor's series for abalone growth option. The method was first introduced by George Adomian who designed the solution to both linear and nonlinear differential equations, including systems of these equations. \cite{[18]}. In addition, The ADM provides an efficient approach for obtaining analytical solutions to a broad and complex set of dynamic systems that represent real-world physical problems. In particular, it offers a well-suited solution for fractional mathematical modelling related to physical issues \cite{[19]}, \cite{[20]}. Also, the method successfully provided in solving the growth model of the density of particles. \cite{[17]}.

The rest of the structure of this study including Sec.2 recalls the preliminaries of the theory. Sec.3 presents the Adomian decomposition method while Sec.4 analyses a new fractional model using the said method, and Sec.5 utilizes the new model to predict the abalone length growth. Finally, Sec.6 discusses the conclusions.

\section{Preliminaries}
This section is a brief review of some key definitions and properties of fractional integrals and derivatives.
\begin{dfn} \cite{[21]}
Let $g:[0,\infty)\longrightarrow\mathbb{R}$ is a continuous function. The definition of the Riemann-Liouville fractional integral and derivative of order $\beta$ are given as.
\begin{equation}
I^\beta_a g(s)=\frac{1}{\Gamma(\beta)}\int^s_a (s-\xi)^{\beta-1}g(\xi)d\xi
\end{equation} 
and
\begin{equation}
^{RL}D^\beta_s g(s)=\begin{cases}\frac{d^n}{ds^n}    \text{                   if    }\beta=m\in\mathbb{N},\\ \\ \frac{d^n}{d s^n}\int^s_0\frac{(s-\xi)^{n-\beta-1}}{\Gamma(n-\beta)}g(\xi)d\xi\text{             if    }m-1<\beta<m, m\in\mathbb{N}.\end{cases}
\end{equation}
\end{dfn}
However, the Equation (7) contradicts with the fundamental theories of derivative due to derivative of the constant function is zero for $\beta\in(m-1,m)$. Therefore, the Definition 2.2 is provided as follows.
\begin{dfn}( \cite{[22]}, \cite{[23]})  
Let a real function $g$ is continuous on $[0,\infty)$ . The derivative of a function $g$ with order $\beta$ is 
\begin{equation}
^{C}D^\beta_a g(s)=\begin{cases}\frac{d^ng(s)}{ds^n}    \text{                   if    }\beta=m\in\mathbb{N}, \\  \\ \int^s_a\frac{(s-\xi)^{m-\beta-1}}{\Gamma(m-\beta)}g^{(m)}(\xi)d\xi \text{             if    }m-1<\beta<m, m\in\mathbb{N}.\end{cases}
\end{equation}
\end{dfn}
Clearly, the Eq.(8) indicates that it does not contradict with fundamental theories of derivative and it is called the Caputo fractional derivative definition. Furthermore, some basic properties are provided as below.
\begin{enumerate}
\item $\left(I^\alpha_a I^\beta_a g\right)(s)=\left(I^\beta_a I^\alpha_a g\right)=\left(I^{\alpha+\beta}g\right)(s)$,
\item $I^\alpha_a(s-a)^\gamma=\frac{\Gamma(\gamma+1)}{\Gamma(\gamma+\alpha+1)}(s-a)^{\gamma+\alpha}$,
\item $\left({I^\alpha_a} ^{C}D^\alpha_a g\right)(s)=g(s)-\sum\limits^{m-1}_{k=0} g^{(k)} (a)\frac{(s-a)^k}{k!}$,
\end{enumerate}
with $\alpha, \beta >0$, $a\geq 0$, $m-1<\alpha<m$,$m\in\mathbb{N}$, and $\gamma>-1$\cite{[24]}.\\

In the same way, the theory can be generalized into fractional partial derivative presented in Definition 2.3.
\begin{dfn}\cite{[25]}
Let $C,D\subseteq [0,\infty)$ and a function $g:C\times D\longrightarrow\mathbb{R}$ is continuous on $C\times D$. The partial derivatives of a function $g$ with order $\beta$ are defined as
\begin{equation}
\frac{\partial^\beta g(s,t)}{\partial s^\beta}=\frac{1}{\Gamma(m-\beta)}\int\limits^s_0 (s-\xi)^{m-\\beta-1}\frac{\partial^m g(\xi,t)}{\partial \xi^m} d\xi, \text{    } \beta\in (m-1, m)
\end{equation}
and
\begin{equation}
\frac{\partial^\beta g(s,t)}{\partial t^\beta}=\frac{1}{\Gamma(m-\beta)}\int\limits^t_0 (t-\xi)^{m-\beta-1}\frac{\partial^m g(s,\xi)}{\partial \xi^m} ds, \text{    } \beta\in (m-1, m).
\end{equation}
\end{dfn}
Thus, the fundamental properties can be similarly generalized. In addition, we provide the definition of the Mittage-Leffler function which involved in solving the fractional differential equations as follows.

\begin{dfn}\cite{[26]}
The Mittage-Leffler function with parameter $\alpha$ is defined by the series expansion
\begin{equation}
E_{\alpha}(z)=\sum\limits^{\infty}_{m=0}\frac{z^m}{\Gamma(m\alpha+1)},  \text{        } \alpha>0, z\in\mathbb{C}
\end{equation}
which the series is convergent. This series is a simple generalization of the exponential function. \end{dfn}
Similarly,
\begin{dfn}\cite{[26]}
The Mittage-Leffler function of two-parameter is defined as: 
\begin{equation}
E_{\alpha, \beta}(z)=\sum\limits^{\infty}_{m=0}\frac{z^m}{\Gamma(m\alpha+\beta)},  \text{        } \alpha, \beta>0, z\in\mathbb{C}
\end{equation}
which the series is convergent and where $\Gamma(\cdot)$ is the Gamma function.
\end{dfn}

\section{Adomian Decomposition Method}
This section discusses about the brief technique of the Adomian decomposition method.  Generally, the technique is converting the differential equations into a linear operator and a nonlinear operator \cite{[17]}, \cite{[27]}. To clarify the ADM technique, let us recall the equation below,
\begin{equation}
\frac{\partial w(s,t)}{\partial t}+\frac{\partial^{\beta} w(s,t)}{\partial s^{\beta}}+w(s,t)+(w(x,t))^2=g(s,t)
\end{equation}
where $n-1<\alpha<n$, $n\in\mathbb{N}$. \\
By ADM, the equation (12) is written in operator form as follows
\begin{equation}
L_t w+L_s w+N w=g(s,t),
\end{equation} 
where the linear operator of $L_t=\frac{\partial}{\partial t}$, and $L_s=\frac{\partial^{\beta}}{\partial s^{\beta}}$ are invertible and $N(w)$ deputizes the non-linear term, that is $N(w)=w(s,t)+(w(s,t))^2$.\\ 
Then the equation (13) is written as
\begin{equation}
w=\varphi+L^{-1}_t g(s,t)-L^{-1}_t(L_s w)-L^{-1}_t (N(w)),
\end{equation}
with $\varphi$ is representing the function arising from initial condition. The nonlinear operator $N(w)$ is usually expressed by an infinite series known as the Adomian polynomials, 
\begin{equation*}
N(w)=\sum\limits^{\infty}_{n=0}A_n,
\end{equation*}
with
\begin{equation}
A_n=\frac{1}{n!}\frac{d}{d\lambda^n}N\left(\sum\limits^n_{i=0}\lambda^i w_i\right)_{\lambda=0}.
\end{equation} 
Thus, the solution can be expressed as
\begin{equation}
w(s,t)=\begin{cases}w_0(s,t)=\varphi,\\w_{n+1}(s,t)=L^{-1}_t g(s,t)-L^{-1}_t\left[L_s(w_n)\right]-L^{-1}_t\left(A_n\right)\end{cases}.
\end{equation}


\section{A Growth Model with Fractional Order}
This part presents a new growth model with fractional order as a modification of the McKendrick equation. Here, we show that the model analysed using ADM provides a particular case of Taylor's series, to predict the growth of a species size. For more detail, given $w(s,t)$ represents the density of population and $\eta$ denotes the growth rate. Since the previous model is of first order, we introduce a new model with a fractional-order $\beta$, where $\beta\in(0,1]$, as follows: 
\begin{equation}
\frac{\partial w(s,t)}{\partial t}+\frac{\partial^{\beta}w(s,t)}{\partial s^{\beta}}=\eta w(s,t),
\end{equation} 
with the initial condition
\begin{equation*}
w(s,0)=Me^{rs},
\end{equation*}
where $M$ is the initial number of individuals and $r$ denoting the initial rate of growth. To apply the ADM, we deliver the Eq.(17) in operator form as follows
\begin{equation}
L_t w+L_s w(s,t)=\eta w
\end{equation} 
where $L$ is a linear operator invertible. By setting,
\begin{equation*}
w=\sum\limits^{\infty}_{n=0}w_n,
\end{equation*}
and $w_0=Me^{rs}$, we obtain 
\begin{equation}
w(s,t)=\begin{cases}w_0(s,t)=Me^{rs},\\w_{n+1}(s,t)=-L^{-1}_t\left[L_s(w_n)\right]+\eta L^{-1}[w_n]\end{cases}.
\end{equation}
Using the Caputo fractional derivative, the running iteration of Eq.(19) provides
\begin{equation}
w_1(s,t)=-L^{-1}_t\left[L_s(w_0)\right]+\eta L^{-1}[w_0]=(\eta-r^{\beta})Me^{rs}t.
\end{equation}
The same manner, we find that
\begin{equation}
w_2(s,t)=(\eta-r^{\beta})^2 Me^{rs}\frac{t^2}{2!},
\end{equation}
\begin{equation}
w_3(s,t)=(\eta-r^{\beta})^3 Me^{rs}\frac{t^3}{3!},
\end{equation}
and the proceeding further the $w_n(s,t)$ is expressed by
\begin{equation}
w_n(s,t)=(\eta-r^{\beta})^n Me^{rs}\frac{t^n}{n!}.
\end{equation}
Therefore, the solution of Eq.(15) is 
\begin{equation}
w(s,t)=\lim_{m\to\infty}\sum\limits^m_{n=0}w_n(s,t)=Me^{rs}e^{(\eta-r^{\beta})t}.
\end{equation}
Then the solution will be utilised to predict the abalone growth in the following section. 

\section{Applications}
In this section we present the predicted values of the growth in the length of abalone utilising the Taylor's series of Eq.(24). To use the series, we determine the growth rate of abalone length based on the data that were also used in \cite{[28]}, \cite{[29]}. We calculate the rate of growth using a formula, $\eta=\Delta h/\Delta t$, where $h$ denotes the abalone length. In addition, the initial growth rate $r=0.04305$ which was obtained from previous study \cite{[28]}. By substituting parameter values on the Eq.(24), we obtain the table 1 as follows.

\begin{table}
\centering

\caption{Approximation of Abalone Length}
\begin{tabular}{|l|c|c|c|c|c|c|r|}
\hline
Month & $\eta$ & $h_{0.5}$ & $h_{0.6}$ & $h_{0.7}$ & $h_{0.8}$ & $h_{0.9}$ & $h_1$ \\
\hline
1  & -       & 0.5322 & 0.5322 & 0.5322 & 0.5322 & 0.5322 & 0.5322\\
2  & 0.4936  & 0.7370 & 0.7794 & 0.8119 & 0.8366 & 0.8550 & 0.8687\\
3  & 0.4724  & 1.3924 & 1.4726 & 1.5341 & 1.5805 & 1.6154 & 1.6413 \\
4  & 0.4521  & 1.9934 & 2.1082 & 2.1962 & 2.2627 & 2.3126 & 2.3496\\
5  & 0.4326  & 2.5435 & 2.6900 & 2.8023 & 2.8872 & 2.9508 & 2.9981\\
6  & 0.4239  & 3.0768 & 3.2540 & 3.3898 & 3.4925 & 3.5694 & 3.6267\\
7  & 0.3962  & 3.5397 & 3.7436 & 3.8998 & 4.0179 & 4.1065 & 4.1723\\
8  & 0.0380  & 3.9611 & 4.1892 & 4.3641 & 4.4963 & 4.5954 & 4.6691\\
9  & 0.3628  & 4.3558 & 4.6067 & 4.7989 & 4.9444 & 5.0533 & 5.1344\\
10 & 0.3472  & 4.7240 & 4.9960 & 5.2045 & 5.3622 & 5.4804 & 5.5683\\
11 & 0.3322  & 5.0643 & 5.3559 & 5.5794 & 5.7485 & 5.8751 & 5.9694\\
12 & 0.3179  & 5.3797 & 5.6895 & 5.9269 & 6.1065 & 6.2410 & 6.3412\\
13 & 0.3043  & 5.6726 & 5.9993 & 6.2497 & 6.4390 & 6.5809 & 6.6865\\
14 & 0.2911  & 5.9436 & 6.2859 & 6.5482 & 6.7467 & 6.8953 & 7.0059\\
15 & 0.2786  & 6.1961 & 6.5529 & 6.8264 & 7.0332 & 7.1882 & 7.3035\\
16 & 0.2666  & 6.4308 & 6.8011 & 7.0849 & 7.2996 & 7.4604 & 7.5801\\
17 & 0.2551  & 6.6491 & 7.0321 & 7.3255 & 7.5475 & 7.7138 & 7.8375\\
18 & 0.2443  & 6.8540 & 7.2487 & 7.5512 & 7.7800 & 7.9515 & 8.0790\\
19 & 0.2336  & 7.0429 & 7.4485 & 7.7593 & 7.9944 & 8.1705 & 8.3016\\
20 & 0.2236  & 7.2206 & 7.6365 & 7.9551 & 8.1962 & 8.3768 & 8.5111\\
21 & 0.2140  & 7.3866 & 7.8120 & 8.1380 & 8.3846 & 8.5693 & 8.7068\\
22 & 0.2047  & 7.5410 & 7.9753 & 8.3081 & 8.5599 & 8.7485 & 8.8888\\
23 & 0.1960  & 7.6869 & 8.1297 & 8.4689 & 8.7255 & 8.9178 & 9.0608\\
24 & 0.1875  & 7.8225 & 8.2730 & 8.6182 & 8.8793 & 9.0750 & 9.2205\\
\hline
\end{tabular}
\end{table}
\newpage
\normalsize
Consider that $h_{\beta}$ represents the abalone length with order $\beta$. By comparing the lengths of abalone on the Table 2 with the real data, we find that the mean absolute error of $h_{0,5}, h_{0.6}, h_{0.7}, h_{0.8}, h_{0.9}$, and $h_1$ are 0.2622; 0.5373; 0.7517; 0.9155; 1.0382; and 1.1294 respectively. Hence, the optimal result is achieved with a fractional order of $\beta=0.5$. It can also be shown in Figure 1 below:

\begin{figure}[h]
\centering
\includegraphics[width = 1.0\textwidth]{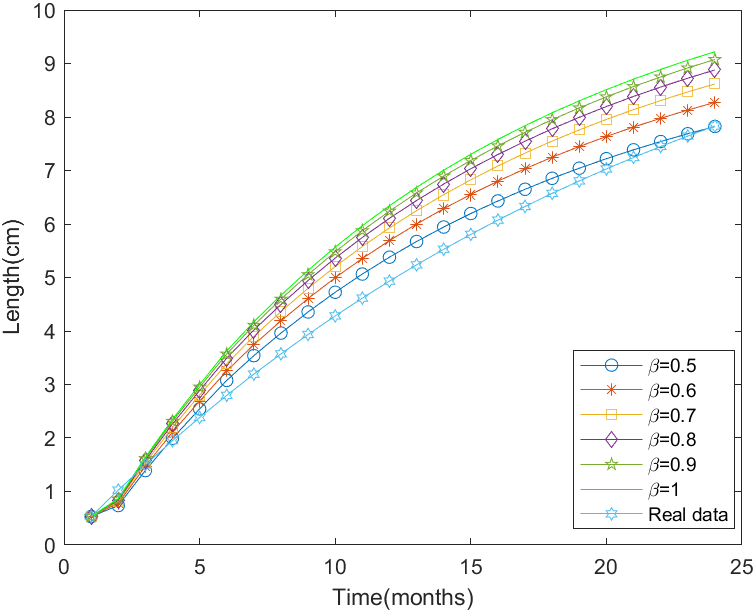}
\caption{Abalone length comparison for 24 months of observation data}
\end{figure}
\newpage

The Figure 1 shows that as the rate of the abalone length growth and the fractional order decrease, the values of the abalone length approach to the real data. The chart with order $\beta=0.5$ is the closest to the real data compared with the other charts. Moreover, note that in the 24th month, the abalone length approximation with $\beta=0.5$ is almost equal to the real data.    

On the other side, the following system equations can also be used to express the Equation (17).
\begin{equation*}
\frac{\partial h(s,t)}{\partial t}+\frac{\partial^{\beta}h(s,t)}{\partial s^{\beta}}=\eta_1 h(s,t),
\end{equation*}
\begin{equation*}
\frac{\partial h(s,t)}{\partial t}+\frac{\partial^{\beta}h(s,t)}{\partial s^{\beta}}=\eta_2 h(s,t)
\end{equation*}
\begin{equation*}
\vdots
\end{equation*}
\begin{equation*}
\frac{\partial h(s,t)}{\partial t}+\frac{\partial^{\beta}h(s,t)}{\partial s^{\beta}}=\eta_{23} h(s,t)
\end{equation*}
or it can be written as
\\
$$
\left[\begin{array}{cc}
1 & 1 \\
1 & 1 \\
\vdots & \vdots \\
1 & 1
\end{array}\right]
\ \ 
\left[\begin{array}{c}
 \displaystyle \frac{\partial h(s,t)}{\partial t} \\
\displaystyle \frac{\partial^{\beta}h(s,t)}{\partial s^{\beta}} \\
 \end{array}\right]
\ = \ \left[\begin{array}{cc}
\eta_1 h(s,t) \\
\eta_2 h(s,t) \\
\vdots \\
\eta_{23} h(s,t)
\end{array}\right]
$$
\\and
\begin{equation*}
f(\eta_1, \beta)= h_1, f(\eta_2, \beta)= h_2, \ldots , f(\eta_{23}, \beta)=h_{23}.  
\end{equation*}

\section{Conclusion}
This research showed that modifying the McKendrick model into a new growth model with fractional order, using an exponential function as the initial condition, can accurately predict the abalone length growth. The new model was analysed using the Adomian decomposition method, resulting to a Taylor series. By substituting the parameter values of this series which were obtained from real data and simulated with various fractional orders, the fractional order $\beta=0.5$ provides results that best reflects the real data, compared to the other orders, including integer orders. For future research, the model can be modified by involving some variables such as competition interaction, diseases, and temperature. \\
\\



\begin{thebibliography}{99}
\bibitem{[1]} Romdhini U.M.\textit{ Optimalization strategy of abalone (haliotis asinina) cultivation in Lombok-island using the matrix Leslie}: INA-Rxiv 2018/3.
\bibitem{[2]} Samson L. M., Shikaa S. \textit{Mathematical modeling of Taraba State population growth using exponential and logistic models}: Results in Control and Optimization 12 (2023) 100265

\bibitem{[3]} Tejera-V. et al., \textit{Estimated Effect of COVID-19 Lockdown on Skin Tumor Size and Survival: An Exponential Growth Model}: Actas Dermosifiliogr. 2020;111(8):629---638.

\bibitem{[4]} Mohsin M., Zaidi A.A., Brunt B. \textit{Dynamics of cell growth: Exponential growth and division after a minimum cell size}: Partial Differential Equations in Applied Mathematics 11 (2024) 100814.

\bibitem{[5]} Kamdoum O.J., Nouadjep S.N., Ndinakie G.P., \textit{Predicting electricity demand in Cameroon with enhanced modified exponential models: A focus on the southern integrated grid (S.I.G)}: Prime - Advances in Electrical Engineering, Electronics and Energy 7 (2024) 100491.

\bibitem{[6]} Kwong C.K. Paul.  \textit{A Time-Dependent Mckendrick Population Model for Logistic Transition}: Mathl. Comput. Modelling Vol. 15, No. 10, pp. 49-59, 1991.

\bibitem{[7]} Gourley A. S., Liu R. \textit{Delay model for populations that experience competition at immature life stages}: J. Differential Equations 259 (2015) 1757-1777. 

\bibitem{[8]} Halder J., Tumuluri K.S., \textit{A higher order numerical scheme to a nonlinear McKendrick–Von Foerster equation with singular mortality}: Applied Numerical Mathematics, volume 202, August 2024, pages 21-41.

\bibitem{[9]} Doumic M., Perthame B., Ribes E., Salort D., Toubiana N.  Toward an integrated workforce planning framework using structured equations: European Journal of Operational Research 262 (2017) 217--230.
 
\bibitem{[10]} Górajski M., Machowska D. \textit{The effects of technological shocks in an optimal goodwill model with a random product life cycle}: Computers and Mathematics with Applications 76 (2018) 905--922.

\bibitem{[11]} Anzai A., Kawatsu L., Uchimura K., Nishiura H. \textit{Reconstructing the population dynamics of foreign residents in Japan to estimate the prevalence of infection with Mycobacterium tuberculosis}: Journal of Theoretical Biology 489 (2020) 110160.

\bibitem{[12]} Adimy M., Chekroun A., Kuniya T.  \textit{Traveling waves of a differential-difference diffusive Kermack-McKendrick epidemic model with age-structured protection phase}: J. Math. Anal. Appl. 505 (2022) 125464.

\bibitem{[13]} Heneghan F. R. et.al.  \textit{A functional size-spectrum model of the global marine ecosystem that resolves zooplankton composition}: Ecological Modelling 435 (2020) 109265.

\bibitem{[14]} Yunfei L., Pei Y., Yuan R. \textit{On a periodic age-structured population model with spatial structure: Non-linear Analysis}: Real World Applications 61 (2021) 103337.

\bibitem{[15]} Veprauskas A. \textit{A non-linear continuous-time model for a semelparous species}: Mathematical Biosciences 297 (2018) 1--11.

\bibitem{[16]} Liu K., Lou Y., Wu J. \textit{ Analysis of an age structured model for tick populations subject to seasonal effects}: J. Differential Equations 263 (2017) 2078--2112.

\bibitem{[17]} Arora G., Hussain S., Kumar R., \textit{Comparison of variational iteration and Adomian decomposition methods to solve growth, aggregation and aggregation-breakage equations}: Journal of Computational Science 67 (2023) 101973.

\bibitem{[18]} Cocom L.B., Estrella G.A., Vales A.E. \textit{Solving delay differential systems with history functions by the Adomian decomposition method}: Applied Mathematics and Computation 218 (2012) 5994–6011.

\bibitem{[19]} Ray S.S., Bera K.R. \textit{An approximate solution of a nonlinear fractional differential equation by Adomian decomposition method}: Applied Mathematics and Computation 167 (2005) 561–571.

\bibitem{[20]} Ray S.S., Bera K.R. \textit{Analytical solution of a fractional diffusion equation by Adomian decomposition method}: Applied Mathematics and Computation 174 (2006) 329–336.

\bibitem{[21]} Ortigueira D.M., Machado T.A.J. \textit{ What is a fractional derivative?}: Journal of Computational Physics 293 (2015) 4--13.

\bibitem{[22]} Tavares D., Almeida R., Torres M.F.D.  \textit{Caputo derivatives of fractional variable order: Numerical approximations}: Commun Nonlinear Sci Numer Simulat 35 (2016) 69–87.

\bibitem{[23]} Sivalingam SM., Govindaraj V. \textit{A novel numerical approach for time-varying impulsive fractional differential equations using theory of functional connections and neural network}: Expert Systems with Applications 238 (2024) 121750.

\bibitem{[24]} Momani S., Odibat Z., Erturk S.V. \textit{Generalized differential transform method for solving a spaceand time-fractional diffusion-wave equation}: Physics Letters A 370 (2007) 379--387.

\bibitem{[25]} Entezari M., Abbasbandy S., Babolian E. \textit{Numerical Solution of Fractional Partial Differential Equations with Normalized Bernstein Wavelet Method}: Appl. Appl. Math. ISSN: 1932--9466; Vol. 14, Issue 2 (December 2019), pp. 890 -- 909.

\bibitem{[26]} Nazir Amna et al. \textit{On generalized fractional integral with multivariate Mittag-Leffler function and its applications}: Alexandria Engineering Journal (2022) 61, 9187–9201. 

\bibitem{[27]}Javed I., Iqbal S., Ali J., Siddique I., Younas M. H. \textit{Unveiling the intricacies: Analytical insights into time and space fractional order inviscid burger’s equations using Adomian decomposition method}: Partial Differential Equations in Applied Mathematics 11 (2024) 100817.

\bibitem{[28]} Susanto M., Romdhini U.M., Kamali R.S., Zurfani L. \textit{Logistic model of abalone’s length growth in Sekotong, West Lombok}: AIP Publishing LLC-2199 (2019)-030002.

\bibitem{[29]} Susanto M., Kilicman A., Nadihah W. \textit{Fractional growth model of abalone length}: Partial Differential Equations in Applied Mathematics 10 (2024) 100668.



\end{thebibliography}
\end{document}